\documentclass[11pt]{amsart}
\usepackage[margin=28mm]{geometry}
\usepackage{amsmath,amssymb,amsthm}
\usepackage[hidelinks]{hyperref}

\newtheorem{theorem}{Theorem}[section]
\newtheorem{proposition}[theorem]{Proposition}
\newtheorem{lemma}[theorem]{Lemma}

\theoremstyle{remark}
\newtheorem{remark}[theorem]{Remark}
\theoremstyle{definition}
\newtheorem{definition}[theorem]{Definition}

\newtheorem*{thmA}{Theorem A}
\newtheorem*{thmB}{Theorem B}
\newtheorem*{corC}{Corollary C}

\newcommand{\R}{\mathbb R}
\newcommand{\Z}{\mathbb Z}
\newcommand{\HH}{\mathbb H}
\newcommand{\Sp}{\mathrm{Sp}}
\newcommand{\SO}{\mathrm{SO}}
\newcommand{\Oo}{\mathrm{O}}
\newcommand{\Ad}{\mathrm{Ad}}
\newcommand{\im}{\operatorname{Im}}
\newcommand{\coker}{\operatorname{coker}}
\newcommand{\sect}{\operatorname{sec}}
\newcommand{\st}{\star}
\newcommand{\Th}{\Theta}
\newcommand{\pis}[1]{\pi_{#1}^{s}}

\title[Positive sectional curvature on all smooth 13-spheres]{Positive sectional curvature on all smooth homotopy 13-spheres}
\author[Z. Chen]{Zhiqi Chen}
\address{School of Mathematics and Statistics,\newline Guangdong University of Technology,\newline Guangzhou 510520, P. R. China}
\email{chenzhiqi@nankai.edu.cn}

\author[S. Deng]{Shaoqiang Deng}
\address{School of Mathematical Sciences and LPMC,\newline Nankai University,\newline Tianjin 300071, P. R. China}
\email{dengsq@nankai.edu.cn}

\author[J. Ni]{Junna Ni}
\address{School of Mathematical Sciences,\newline South China Normal University,\newline Guangzhou 510631, P. R. China}
\email{nijunna@126.com}

\hypersetup{pdftitle={Positive sectional curvature on all smooth homotopy 13-spheres},pdfauthor={Zhiqi Chen, Shaoqiang Deng, Junna Ni}}

\date{November 7, 2026}

\begin{document}

\begin{abstract}
Every smooth homotopy $13$-sphere admits a Riemannian metric of strictly positive sectional curvature. We construct a special $S^3$--$S^3$ bundle over $S^{13}$ whose star quotient represents the nonzero class $\pm\alpha_1\beta_1\in\Theta_{13}\cong\mathbb Z/3$. The defining equivariant map factors through the Blakers--Massey map by a smooth Hopf construction representing an unstable $J$-generator. This factorisation retains the prescribed base representation and makes Schultz's plumbing composition formula applicable; a normal-invariant computation in smoothing theory then identifies the smooth structure. On the quotient we implement the compatible-disk strategy with complete curvature estimates on arbitrary two-planes, matching boundary metrics, and a strictly positive sum of outward second fundamental forms. A direct smoothing argument completes the metric construction. Orientation reversal gives the remaining nonzero class in $\Theta_{13}$.
\end{abstract}

\maketitle

\section{Introduction}\label{sec:introduction}

We prove that every smooth homotopy $13$-sphere admits a metric of strictly positive sectional curvature. The construction joins an explicit equivariant model of a nonzero class in $\Th_{13}$ to a compatible pair of positively curved quotient disks. Its topological step is a factorisation through an unstable $J$-image which preserves the representation required by the geometric model.

Recent preprints construct positive sectional curvature on all homotopy $7$-spheres \cite{HLY}, on the exotic $8$-sphere and the order-three $10$-spheres \cite{GG}, and on an even subgroup of the $11$-spheres \cite{DHZ}. Galaz-Garc\'{\i}a constructs exotic positively curved spheres in infinitely many dimensions congruent to $3$ modulo $4$ \cite{GG2}; that family does not include dimension $13$. The present construction treats this dimension by combining Speran\c{c}a's bundle models \cite{Sp} with a Hopf factorisation adapted to the source action.

Kervaire and Milnor's calculation $\Th_{13}\cong\Z/3$ \cite{KM} reduces the problem to constructing a metric on one nonstandard sphere. We separate the geometric assertion from the identification of its smooth structure. A special $S^3$--$S^3$ bundle is a principal $S^3$-bundle with a second commuting free $S^3$-action, denoted by $\star$; precise conventions are fixed in Section~\ref{sec:preliminaries}.
\begin{thmA}
There is a special $S^3$--$S^3$ bundle $E^{16}\to S^{13}$, obtained by a smooth equivariant pullback of the Gromoll--Meyer bundle $\Sp(2)\to S^7$, such that its star quotient
\[
\Sigma^{13}=E^{16}/S^3_\star
\]
admits a smooth Riemannian metric of strictly positive sectional curvature.
\end{thmA}

The same quotient has a nonzero smooth-structure invariant. Here $\alpha_1$ and $\beta_1$ are the order-three stable classes in stems $3$ and $10$, with conventions recalled in Section~\ref{sec:preliminaries}.
\begin{thmB}
The quotient $\Sigma^{13}$ in Theorem~A is an exotic homotopy sphere. Under the Kervaire--Milnor identification
\[
\Th_{13}\cong\pis{13}\cong\Z/3,
\]
its class is $\pm\alpha_1\beta_1$.
\end{thmB}

The two nonzero elements of $\Th_{13}$ differ by orientation reversal, which preserves sectional curvature.
\begin{corC}
Every smooth homotopy $13$-sphere admits a smooth Riemannian metric of strictly positive sectional curvature.
\end{corC}
\begin{proof}
The standard sphere has its round metric. By Theorem~B the other two oriented classes are $[\Sigma^{13}]$ and $[-\Sigma^{13}]$. Theorem~A provides a metric on $\Sigma^{13}$, and the same metric on the underlying smooth manifold has the same sectional curvature after its orientation is reversed.
\end{proof}

\subsection*{The equivariant model and its smooth structure}
The construction uses
\[
U=\R^6\oplus\HH_L,
\]
with trivial action on $\R^6$ and left multiplication on $\HH_L$. The required map $a_{13}:S(U)\to S^3$ must be conjugation-equivariant and must produce a nonzero plumbing class. We obtain both properties from the actual factorisation
\[
a_{13}=b\circ g,\qquad
g:S(\R^6_{\mathrm{triv}}\oplus\HH_L)\longrightarrow
S(\R^3_{\mathrm{triv}}\oplus\HH_L),\qquad
[g]=J([\R^2\oplus L])=\nu_6.
\]
Here $b:S^6\to S^3$ is the Blakers--Massey map and $g$ is a smoothed Hopf construction. The domain and target actions are compatible with this formula, and the radial profiles give smooth maps at both join ends. Retaining the representation $U$ keeps the second reduced clutching class equal to a generator; changing that representation can annihilate the plumbing class, as explained in Remark~\ref{rem:choice}.

Speran\c{c}a's plumbing description \cite{Sp} and Schultz's composition formula \cite{Sch} then give
\[
\Sigma^{13}=\pm\tau_{3,10}(\gamma,\Sigma^{10}),\qquad
\gamma=[\R^2\oplus L],
\]
where $\Sigma^{10}$ is Speran\c{c}a's order-three sphere. The Samelson obstruction to applying the composition formula vanishes in $\pi_6(\SO(7))$. Schultz's description of the Milnor--Munkres--Novikov pairing in $\pi_*(\mathrm{PL}/\Oo)$ \cite{Sch}, together with the identification of the Kervaire--Milnor invariant as a normal invariant \cite{Brumfiel}, identifies the class as $\pm\nu\beta_1=\pm\alpha_1\beta_1$. Thus the same explicit map supplies both the equivariance needed for the bundle and the homotopy class detecting exoticity.

\subsection*{The metric construction and its relation to earlier work}
We implement the compatible-disk strategy of \cite{HLY,GG}. The general transfer statement used in this setting is already present in \cite[Theorem~4.1]{GG}. Our geometric part gives a complete proof for the bundle above, including the estimates and smoothing needed to pass from its two quotient disks to a metric on the prescribed smooth sphere.

On the southern bundle a variable fibre radius supplies positive mixed curvature, and a decomposable-bivector calculation controls the connection terms on arbitrary two-planes. On the northern bundle a doubly warped metric is estimated on every plane perpendicular to the star orbits. The source boundary metrics agree in the bundle trivialisation, so the quotient boundary metrics match under the actual attaching map. Their outward second fundamental forms have strictly positive sum. A direct strict gluing argument then produces a smooth metric with positive sectional curvature. All curvature calculations and parameter choices are included in Section~\ref{sec:geometry}; the metric proof uses no geometric existence theorem from the recent preprints.

Section~\ref{sec:preliminaries} fixes notation and the classical inputs. Section~\ref{sec:construction} constructs the equivariant bundle, Section~\ref{sec:topology} proves Theorem~B, and Section~\ref{sec:geometry} proves Theorem~A.

\section{Preliminaries}\label{sec:preliminaries}

We fix the topological inputs, representation conventions and curvature signs used below.

\subsection{Homotopy spheres and stable classes}
For $n\ge5$, let $\Th_n$ be the group of oriented diffeomorphism classes of smooth homotopy $n$-spheres under connected sum. Its zero element is the standard sphere and its inverse is orientation reversal. Write
\[
\pis j=\mathop{\rm colim}_{m}\pi_{m+j}(S^m),\qquad
\coker J_j=\pis j/\operatorname{im}J_j.
\]
We use stable composition $\pis i\otimes\pis j\to\pi_{i+j}^s$ and the unstable and stable $J$-homomorphisms
\[
J:\pi_r(\SO(p))\longrightarrow\pi_{p+r}(S^p),\qquad
J:\pi_r(\SO)\longrightarrow\pis r.
\]
Here $\SO$ is the stable orthogonal group.

The Kervaire--Milnor homomorphism $\chi_n:\Th_n\to\coker J_n$ has kernel $bP_{n+1}$, the subgroup of spheres bounding parallelisable manifolds \cite{KM}. The low-dimensional data needed here are
\[
\begin{array}{c|c|c|c|c}
n&\Th_n&\pi_n^s&bP_{n+1}&\pi_n(\SO)\\\hline
10&\Z/6&\Z/6&0&0\\
13&\Z/3&\Z/3&0&0.
\end{array}
\]
For the stable stems and orthogonal groups see \cite[\S1.1, pp.~3--5]{Ra}; see also \cite{Toda}. Thus $\chi_{10}$ and $\chi_{13}$ are isomorphisms onto the corresponding stable stems: their kernels and the images of $J$ vanish, and source and target have the same order.

Choose a generator $\nu\in\pis3\cong\Z/24$ by the $J$-normalisation in Lemma~\ref{lem:a13}. Denote its $3$-primary component by $\alpha_1$, and choose a generator $\beta_1$ of the order-three subgroup of $\pis{10}$. The product
\[
\alpha_1\beta_1\ne0\quad\hbox{generates }\pis{13}\cong\Z/3.
\]
This is the case $p=3$, $q=2p-2=4$ of \cite[Theorem~1.1.14, p.~5]{Ra}; see also \cite[Figure~1.2.19, p.~13]{Ra}.

\subsection{Quaternionic actions and the Blakers--Massey map}
Identify $S^3$ with the unit quaternions. We write
\[
L(q)w=qw,\qquad \mathrm{conj}(q)y=qy\bar q,
\qquad \mathrm{conj}=\R\oplus\Ad.
\]
The notation $\HH_L$ denotes the left representation, and $\R^k$ always denotes a trivial representation. For $m\ge5$ we normalise $\pi_3(\SO(m))\cong\Z$ by $[L]=1$; then $[\Ad]=2$.

On $S^6=S(\im\HH\oplus\HH)$ the Blakers--Massey map is
\begin{equation}\label{eq:BM-preliminary}
b(p,w)=
\begin{cases}
\dfrac{w}{|w|}\exp(\pi p)\dfrac{\bar w}{|w|},&w\ne0,\\[4pt]
-1,&w=0.
\end{cases}
\end{equation}
It represents a generator of $\pi_6(S^3)\cong\Z/12$ \cite[\S2]{Sp} and satisfies
\[
b(p,qw)=q\,b(p,w)\,q^{-1},\qquad b(0,w)=1\quad(|w|=1).
\]
For smoothness at $w=0$, put $s=|p|^2=1-|w|^2$ and write
\[
b(p,w)=\cos(\pi\sqrt s)+\psi(s)wp\bar w,\qquad
\psi(s)=\frac{\sin(\pi\sqrt s)}{\sqrt s(1-s)}.
\]
The function $\psi$ extends smoothly to $[0,1]$, with endpoint values $\pi$ and $\pi/2$. This expression proves smoothness on the whole sphere.

For based maps $\alpha:S^p\to G$, $\beta:S^q\to G$ into a topological group, the Samelson product $\langle\alpha,\beta\rangle\in\pi_{p+q}(G)$ is represented by
\[
(x,y)\longmapsto\alpha(x)\beta(y)\alpha(x)^{-1}\beta(y)^{-1}
\]
on $S^p\wedge S^q\cong S^{p+q}$. We use it for orthogonal groups in Section~\ref{sec:topology}.

\subsection{Equivariant bundles and joins}
A map is equivariant if $f(q\cdot x)=q\cdot f(x)$. For a map to $S^3$ with conjugation action this means $f(q\cdot x)=qf(x)q^{-1}$. We use the sphere model
\[
S^k*S^l=S(\R^{k+1}\oplus\R^{l+1}),
\qquad (a*d)[x,y,t]=[a(x),d(y),t].
\]
An equivariant join is continuous; the flat radial formula in Section~\ref{sec:construction} supplies the smooth representative used for pullback.

For a principal bundle $\operatorname{pr}:P\to Y$, its pullback by $f:X\to Y$ is
\[
f^*P=\{(x,Q):f(x)=\operatorname{pr}(Q)\}.
\]
If a second action on $P$ covers the action on $Y$ and commutes with the principal action, equivariance of $f$ makes $(x,Q)\mapsto(q\cdot x,q\star Q)$ an action on $f^*P$ with the same commutation property.

\begin{definition}\label{def:special-star}
A special $S^3$--$S^3$ bundle is a principal right $S^3$-bundle equipped with a commuting free left $S^3$-action, called the star action. When its base is $S^n=S(\R e\oplus V)$ and the induced star action on the base is the restriction of an orthogonal representation fixing $e$, we call it a star bundle.
\end{definition}

\subsection{Curvature and boundary conventions}\label{subsec:curvature-conventions}
We use
\[
R(A,B)C=\nabla_A\nabla_BC-\nabla_B\nabla_AC-\nabla_{[A,B]}C,
\qquad
\sect(A\wedge B)=\frac{\langle R(A,B)B,A\rangle}{|A|^2|B|^2-\langle A,B\rangle^2}.
\]
Write $\mathcal R(\xi,\xi)$ for the numerator when $\xi=A\wedge B$ is decomposable. Our estimates concern actual two-planes; they do not assert positivity of the full curvature operator.

The principal-horizontal distribution is $\ker\omega$ for a principal connection $\omega$. The star-horizontal distribution is the metric orthogonal complement of the star orbits. These are distinct splittings. For orthonormal star-horizontal vectors, O'Neill's formula \cite{ON} is
\[
\sect_{\rm quotient}(d\pi A\wedge d\pi B)
=\sect_{\rm source}(A\wedge B)+3|\mathcal A_AB|^2.
\]
Thus positivity on every star-horizontal source plane implies positivity on the quotient. Section~\ref{sec:geometry} proves positivity on all source planes in the south and on all star-horizontal planes in the north.

The metric $h_{S^k}$ is unit round. With outward unit normal $\nu$, the boundary convention is $B(X,Y)=\langle\nabla_X\nu,Y\rangle$; a small round ball has positive definite $B$.

\section{The equivariant construction}\label{sec:construction}

In this section we construct the bundle in Theorem~A. The main point is to produce an equivariant map with the required nonzero homotopy class while keeping the base representation simple. We first recall the Gromoll--Meyer two-action model and its plumbing interpretation. We then construct the new map and verify its smoothness, equivariance and homotopy class separately. These verifications provide the link between the metric construction and the topological calculation in the next section.

\subsection{The Gromoll--Meyer bundle and Speran\c{c}a's pullbacks}
The Gromoll--Meyer bundle is useful because it already has the two commuting actions required in Definition~\ref{def:special-star}. A point of $\Sp(2)$ is a quaternionic unitary matrix
\[
Q=\begin{pmatrix}a&c\\b&d\end{pmatrix},
\]
and projection to its first column gives $\Sp(2)\to S^7$. The principal right action multiplies the second column by a unit quaternion. The star action is
\[
q\star Q=\begin{pmatrix}qa\bar q&qc\\qb\bar q&qd\end{pmatrix}.
\]
The actions commute by associativity of quaternionic multiplication. The star action is free: if it fixes $Q$, then $qc=c$ and $qd=d$, and at least one of $c,d$ is nonzero. Thus $q=1$. The quotient is the Gromoll--Meyer sphere, and the induced base action is simultaneous conjugation on the two quaternionic coordinates \cite[\S2]{Sp}.

For a join pullback it is convenient to decompose $S^7$ as $D^4\times S^3$ joined to $S^3\times D^4$ along $S^3\times S^3$. In the compatible trivialisation of \cite[\S3]{Sp}, the transition function on this common boundary is $\theta(x,y)=x\bar y$. Let $S^3$ act orthogonally on $S^k\subset\R^{k+1}$ and $S^l\subset\R^{l+1}$ through homomorphisms $\Delta_1$ and $\Delta_2$, respectively, and let $a:S^k\to S^3$ and $d:S^l\to S^3$ be conjugation-equivariant. Denote the pullback by a smooth representative of $a*d$ by $E=(a*d)^*\Sp(2)$. Its transition function is
\begin{equation}\label{eq:r}
r(x,y)=\theta\circ(a\ast d)(x,y)=a(x)\,d(y)^{-1}.
\end{equation}

The transition function is explicit, but we still need a way to read the smooth structure of the star quotient from it. Plumbing provides such a description. Its construction is local: in trivial product neighbourhoods of the zero sections of two disk bundles, one exchanges the base and fibre coordinates and identifies the two products. In the present situation one disk bundle has base $S^{k+1}$ and fibre $D^{l+1}$, while the other has base $S^{l+1}$ and fibre $D^{k+1}$. Both total spaces have dimension $k+l+2$, so the boundary of their plumbing has the same dimension as the star quotient.

The bundles being plumbed are specified by clutching maps. A clutching map records the change of a vector-bundle frame between the two hemispheres of its base sphere. Thus a map $S^k\to\SO(l+1)$ specifies the first disk bundle, and a map $S^l\to\SO(k+1)$ specifies the second. The following result says that the two clutching maps can be read directly from the two representation actions in the join construction. This is the reason for keeping track of those actions, rather than only the ordinary homotopy class of the pullback map.
\begin{theorem}[Speran\c{c}a {\cite[Theorem~3, Corollary~5]{Sp}}]\label{thm:sp}
With $r$ as in \eqref{eq:r}, the star quotient is the twisted double
\[
E/S^3_\st\;\cong\;D^{k+1}\times S^l\;\cup_{\hat r}\;S^k\times D^{l+1},\qquad
\hat r(x,y)=\bigl(r(x,y)\cdot x,\; r(x,y)\cdot y\bigr),
\]
where $r(x,y)\in S^3$ acts on $x$ and $y$ through the given linear actions. It is the boundary of the plumbing of the disk bundles over $S^{k+1}$ and $S^{l+1}$ with clutching maps $\Delta_2 a\in\pi_k(\SO(l+1))$ and $\Delta_1 d\in\pi_l(\SO(k+1))$, where $\Delta_1\colon S^3\to\SO(k+1)$ and $\Delta_2\colon S^3\to\SO(l+1)$ are the homomorphisms defining the two actions.
\end{theorem}

We will use this theorem in two dimensions. The $13$-dimensional quotient is the manifold we seek, while the $10$-dimensional quotient provides a known nonzero class to which the computation can be reduced. It is useful to recall the lower-dimensional example now, so that the later application of Schultz's formula has a concrete topological starting point.

In all join pullbacks below, the join is understood through a smooth equivariant representative, as in \cite[\S3]{Sp}. Speran\c{c}a's exotic $10$-sphere arises in this way: let $b\colon S^6=S(\R^3\oplus\HH)\to S^3$ be the smooth Blakers--Massey map described in \cite[\S2]{Sp}, which satisfies $b(p,qw)=q\,b(p,w)\,\bar q$ for the action $(p,w)\mapsto(p,qw)$ and represents a generator $c=[b]$ of $\pi_6(S^3)\cong\Z/12$. With $S^{10}=S^6\ast S^3$, the action $\rho_{10}(q)(p,w,y)=(p,qw,qy\bar q)$ and $f_{10}=b\ast\mathrm{id}$, put
\[
E^{13}=f_{10}^\ast\Sp(2),\qquad \Sigma^{10}=E^{13}/S^3_\st .
\]
The point of this example is its smooth structure. Its metric will not be used in the argument. The published topological calculation of Speran\c{c}a gives exactly the order-three class needed for the dimension-raising composition formula.
\begin{theorem}[Speran\c{c}a {\cite[Theorem~1]{Sp}}]\label{thm:sp10}
$\Sigma^{10}$ represents a generator of the subgroup of order three of $\Th_{10}\cong\Z/6$.
\end{theorem}

\subsection{A smooth equivariant Hopf construction in dimension nine}
We need a map on $S(\R^6\oplus\HH_L)$ which converts left multiplication on the last coordinate into conjugation on the value. The Blakers--Massey map already has the corresponding property on $S^6$. It is therefore natural to look for an equivariant map into $S^6$ and then compose with $b$. The Hopf construction of an orthogonal representation gives such a map and simultaneously records its ordinary homotopy class as a $J$-image.

To specify the actions, let
\[
\gamma_0:S^3\longrightarrow\SO(6),\qquad
\gamma_0=\R^2\oplus L
\]
act on $\R^6=\R^2\oplus\HH$. This class generates $\pi_3(\SO(6))\cong\Z$, in the normalisation of Section~\ref{sec:preliminaries}. Identify $\R\oplus\R^2$ with $\im\HH$. The target of the desired map is then
\[
S^6=S(\R\oplus\R^6)=S(\im\HH\oplus\HH),
\qquad q\cdot(\lambda,v)=(\lambda,\gamma_0(q)v).
\]
This is exactly the action $(p,v_{\HH})\mapsto(p,qv_{\HH})$ for which $b$ is equivariant. On the source $U=\R^6\oplus\HH$, however, all six coordinates of the first summand are fixed: $q\cdot(z,w)=(z,qw)$. Thus $\gamma_0$ describes the intermediate target action, whereas $\R^6_{\mathrm{triv}}\oplus L$ describes the source action. Keeping these two representations distinct is essential.

The Hopf formula uses the direction of both $z$ and $w$. These directions cease to be defined at the join ends. We eliminate this smoothness problem by making the formula constant on neighbourhoods of those ends. Choose a smooth nonincreasing function $\zeta:[0,1]\to[0,\pi]$, equal to $\pi$ near $0$ and $0$ near $1$, and set
\begin{equation}\label{eq:hopf-map}
g(z,w)=\left(\cos\zeta(|w|),\;
\sin\zeta(|w|)\,\gamma_0\left(\frac w{|w|}\right)\frac z{|z|}\right)
\in S(\R\oplus\R^6).
\end{equation}
For $w=0$ and $z=0$ define the values to be $(-1,0)$ and $(1,0)$, respectively. Those values already hold on whole neighbourhoods, by the choice of $\zeta$.

We now verify both properties needed later: strict equivariance for the bundle pullback, and an ordinary $J$-class for the plumbing calculation. These are properties of the same explicit map.
\begin{lemma}\label{lem:a13}
The map $g:S^9\to S^6$ in \eqref{eq:hopf-map} is smooth and equivariant for the source and target actions just specified. Choose the generator $\nu_6\in\pi_9(S^6)\cong\Z/24$ so that
\begin{equation}\label{eq:g}
[g]=J([\gamma_0])=\nu_6.
\end{equation}
Then the smooth map
\begin{equation}\label{eq:a13}
a_{13}=b\circ g:S(\R^6_{\mathrm{triv}}\oplus\HH_L)\longrightarrow S^3
\end{equation}
satisfies $a_{13}(z,qw)=q\,a_{13}(z,w)q^{-1}$.
\end{lemma}
\begin{proof}
Smoothness holds where both $z,w$ are nonzero by the displayed formula. It holds at the join ends because the map is constant on their neighbourhoods. For equivariance use that $\gamma_0$ is a homomorphism:
\[
\gamma_0\left(\frac{qw}{|w|}\right)\frac z{|z|}
=\gamma_0(q)\gamma_0\left(\frac w{|w|}\right)\frac z{|z|}.
\]
The first coordinate and the profile are unchanged, so this is exactly the desired target action. Equivariance of $b$ then proves the assertion about $a_{13}$.

For the homotopy class, write $w=\cos(\tau/2)u$, $z=\sin(\tau/2)v$ for $0\le\tau\le\pi$, $u\in S^3$, $v\in S^5$. Deform the profile $\zeta(t)$ to $2\arccos t$ by linear interpolation. The endpoint values stay $\pi$ at $t=0$ and $0$ at $t=1$, so this is a continuous equivariant homotopy on the whole sphere, including both ends. Its final map is
\[
[u,v,\tau]\longmapsto
\bigl(\cos\tau,\sin\tau\,\gamma_0(u)v\bigr).
\]
This is the Hopf construction defining $J([\gamma_0])$. The group $\pi_9(S^6)$ is stable by the suspension theorem, and the stable $J$-image of a generator of $\pi_3(\SO)$ generates $\pi_3^s\cong\Z/24$ \cite{Ad}; see also \cite[Theorem~1.1.13]{Ra}. Since $\pi_3(\SO(6))\to\pi_3(\SO)$ is an isomorphism, $J([\gamma_0])$ is a generator as asserted.
\end{proof}

The nonvanishing relevant to exoticity will be proved by the plumbing calculation in Section~\ref{sec:topology}. It is not necessary to compute $[a_{13}]\in\pi_9(S^3)$ separately before constructing the bundle. This order keeps the topological argument focused on the invariant which distinguishes the smooth structures.

\subsection{The 13-dimensional bundle}
The first join factor now has the required equivariant map. For the second factor we use the identity of $S^3$, with conjugation on both its domain and target. Their join therefore maps to the two conjugation factors of the Gromoll--Meyer base.

Let
\[
S^{13}=S(U\oplus\HH),\qquad U=\R^6\oplus\HH,
\]
with the representation
\begin{equation}\label{eq:rho13}
\rho_{13}(q)(z,w,y)=(z,qw,qy\bar q).
\end{equation}
The join $a_{13}\ast\mathrm{id}$ is equivariant. To use an ordinary smooth pullback bundle, we must also address the two join ends. The issue is that the angular variable $x/|x|$ has derivatives containing negative powers of $|x|$. Multiplying by a function which vanishes to every order removes all such powers. This explains the flat function in the next formula. To specify a smooth representative at its join ends, put $h(s)=e^{-1/s^2}$ for $s>0$ and $h(0)=0$. For $x=(z,w)\in U$, $y\in\HH$, $|x|^2+|y|^2=1$, set
\begin{equation}\label{eq:f13}
f_{13}(x,y)=\frac{\bigl(h(|x|)a_{13}(x/|x|),\ h(|y|)y/|y|\bigr)}{\sqrt{h(|x|)^2+h(|y|)^2}}.
\end{equation}
A component is defined to be zero when its corresponding input is zero. Every derivative of the angular factor grows at most as a fixed negative power of the corresponding radius, because $a_{13}$ is smooth on the compact unit sphere. The flatness of $h$ therefore makes the components smooth at the ends, and the denominator never vanishes. Thus $f_{13}$ is smooth and equivariant. Interpolating the positive radial weights $h(s)$ with $s$ gives an equivariant homotopy to the continuous join. Define
\[
E^{16}=f_{13}^*\Sp(2),\qquad \Sigma^{13}=E^{16}/S^3_\st.
\]
The function $h$ is strictly increasing. Hence the preimages of the two Gromoll--Meyer trivialising pieces are the two join pieces separated by $|x|=|y|=1/\sqrt2$. On this separating product of spheres the two weights agree, so the transition map is exactly $r(\xi,\eta)=a_{13}(\xi)\eta^{-1}$, with $(\xi,\eta)\in S^9\times S^3$. Thus Theorem~\ref{thm:sp} applies directly to the smooth map \eqref{eq:f13}. It identifies $\Sigma^{13}$ with the boundary of the plumbing of the $D^4$-bundle over $S^{10}$ with clutching $\mathrm{conj}\circ a_{13}$ and the $D^{10}$-bundle over $S^4$ with clutching $\rho_U=\R^6\oplus L$. Both clutching maps fix a vector, so the reduced clutching maps belong to $\pi_9(\SO(3))$ and $\pi_3(\SO(9))$. This reduction is important: it places both clutching maps in exactly the groups used in the usual homotopy-sphere plumbing pairing. In the next section we compute that pairing and prove that this quotient is exotic.

\section{The diffeomorphism type}\label{sec:topology}

The preceding construction gives a well-defined star quotient. We now identify its smooth structure. The reason for doing this before constructing the metric is that the curvature calculation itself cannot distinguish the standard sphere from an exotic one. The plumbing description makes this distinction accessible to homotopy theory.

There are three steps. We first express the two quotients in a common plumbing notation. We then use Schultz's formula to relate the $13$-dimensional plumbing to the $10$-dimensional example. Finally, we evaluate the resulting class by multiplication in the stable stems.

\subsection{Plumbing pairings}
A plumbing boundary can be recorded as a commutator of two elementary diffeomorphisms on a product of spheres. One diffeomorphism rotates the second factor by an amount depending on the first, and the other rotates the first factor by an amount depending on the second. Their failure to commute records the smooth sphere produced by the plumbing. This description explains both the bilinearity of the pairing and the commutation condition in Schultz's composition formula.

To describe stabilisation, write $\mathbf{1}$ for the identity on an additional real line. Thus $u\oplus\mathbf{1}$ means that the orthogonal map $u$ is extended to act trivially on this line. Following Schultz \cite[\S2]{Sch}, for $u\in\pi_p(\SO(q))$ and $v\in\pi_q(\SO(p))$ let $f_u$, $g_v$ be the diffeomorphisms of $S^p\times S^q$ given by $f_u(x,y)=(x,u(x)y)$ and $g_v(x,y)=(v(y)x,y)$, where $\SO(q)\subset\SO(q+1)$ and $\SO(p)\subset\SO(p+1)$ act on $S^q$ and $S^p$ in the standard way. Choose representatives constant near the appropriate hemispheres, as in \cite[pp.~741--742]{Sch}. The commutator $[f_u,g_v]$ is then supported in a disk and defines the \emph{plumbing pairing}
\[
\sigma_{p,q}\colon\pi_p(\SO(q))\times\pi_q(\SO(p))\longrightarrow\Th_{p+q+1},
\]
which is bilinear \cite[pp.~741--742]{Sch}. Equivalently, $\sigma_{p,q}(u,v)$ is the twisted double $D^{p+1}\times S^q\cup_{[f_u,g_v]}S^p\times D^{q+1}$ \cite[p.~742, Remark]{Sch}; this is the boundary of the plumbing of the disk bundles over $S^{p+1}$ and $S^{q+1}$ with clutching maps $u\oplus\mathbf{1}$ and $v\oplus\mathbf{1}$, which is the form in which Speran\c{c}a states his result \cite[Proposition~4, Corollary~5]{Sp}. In particular this boundary is a homotopy sphere, as recorded by the target $\Th_{p+q+1}$. A second operation is needed to compare dimensions. The Milnor--Munkres--Novikov pairing combines a family of orthogonal rotations with a homotopy sphere, increasing the dimension of the sphere by the dimension of the family. It is the operation that appears when a clutching map is precomposed with a $J$-image. For $\gamma\in\pi_r(\SO)$ and $\Sigma\in\Th_n$ let $\tau_{r,n}(\gamma,\Sigma)\in\Th_{n+r}$ denote this pairing. When $\gamma$ is written in a finite orthogonal group, we use its image under stabilisation; in the application $\pi_3(\SO(6))\to\pi_3(\SO)$ is an isomorphism.

The class of the first clutching map in our $13$-dimensional model is a composite $\alpha\circ g$. Since $g$ is a $J$-image, the following theorem converts that composite into a Milnor--Munkres--Novikov operation on the lower-dimensional plumbing. The condition concerning the other clutching map ensures that the two relevant rotations can be interchanged in this calculation.
\begin{theorem}[Schultz {\cite[Theorem~2.5]{Sch}}]\label{thm:schultz}
Let $\gamma\in\pi_r(\SO(p))$ with $r\le p-2$, $\alpha\in\pi_p(\SO(q))$ and $\beta\in\pi_q(\SO(p))$. Suppose that $\beta$ lies in the image of $\pi_q(\SO(p-r))\to\pi_q(\SO(p))$. Then
\[
\tau_{r,p+q+1}\bigl(\gamma,\sigma_{p,q}(\alpha,\beta)\bigr)=\sigma_{p+r,q}\bigl(\alpha\circ J(\gamma),\beta\bigr),
\]
where $J(\gamma)\in\pi_{p+r}(S^p)$ is the $J$-homomorphism and $\beta$ is regarded in $\pi_q(\SO(p+r))$ via the standard inclusion.
\end{theorem}

Schultz adds that the condition on $\beta$ may be weakened to the vanishing of the Samelson product of $\gamma$ and $\beta$ in $\SO(p+1)$ \cite[p.~746]{Sch}. In his proof the hypothesis is used only to ensure that the diffeomorphisms $(x,y,z)\mapsto(\beta(y)x,y,z)$ and $(x,y,z)\mapsto(\gamma(z)x,y,z)$ of $S^p\times S^q\times S^r$ commute in the group of pseudo-isotopy classes, where $\SO(p)\subset\SO(p+1)$ acts on $S^p$. At the level of homotopy classes, $\gamma$ can be represented in $\SO(r+1)$, and $\SO(p-r)$ and $\SO(r+1)$ are commuting block subgroups of $\SO(p+1)$. Thus under the original hypothesis the induced classes of diffeomorphisms commute.

For our classes, the smaller-group condition in Schultz's original statement will not hold. We therefore use the Samelson-product variant which he explicitly records; an alternative that uses only the original hypothesis is given after \eqref{eq:tau}. The next lemma explains why that variant is natural: a Samelson product is precisely the obstruction to deforming the commutator of these two fibrewise rotations to the identity.
\begin{lemma}\label{lem:samelson}
Let $\beta\in\pi_q(\SO(p+1))$ and $\gamma\in\pi_r(\SO(p+1))$, and let $f_\beta(x,y,z)=(\beta(y)x,y,z)$ and $f_\gamma(x,y,z)=(\gamma(z)x,y,z)$ be the induced diffeomorphisms of $S^p\times S^q\times S^r$. If the Samelson product $\langle\beta,\gamma\rangle\in\pi_{q+r}(\SO(p+1))$ vanishes, then the commutator $[f_\beta,f_\gamma]$ is isotopic to the identity. In particular Theorem~\ref{thm:schultz} holds whenever $\langle\beta,\gamma\rangle=0$ in $\pi_{q+r}(\SO(p+1))$.
\end{lemma}

\begin{proof}
A direct computation gives $[f_\beta,f_\gamma](x,y,z)=\bigl(c(y,z)\,x,\;y,\;z\bigr)$ with
\[
c(y,z)=\beta(y)\gamma(z)\beta(y)^{-1}\gamma(z)^{-1}\in\SO(p+1).
\]
Choosing representatives with $\beta(y_0)=\gamma(z_0)=1$ for base points $y_0,z_0$, the map $c\colon S^q\times S^r\to\SO(p+1)$ is constant equal to $1$ on $S^q\vee S^r$ and therefore factors through the quotient $S^q\times S^r\to S^q\wedge S^r=S^{q+r}$; the resulting based map $\bar c\colon S^{q+r}\to\SO(p+1)$ represents the Samelson product $\langle\beta,\gamma\rangle$ by definition. If $\langle\beta,\gamma\rangle=0$, a based null-homotopy $\bar c_t$ of $\bar c$ composed with the quotient map gives a homotopy $c_t\colon S^q\times S^r\to\SO(p+1)$ from $c$ to the constant map $1$. By relative smooth approximation on $S^q\times S^r\times[0,1]$, the homotopy $c_t$ may be taken smooth with its endpoint maps fixed. Then $F_t(x,y,z)=(c_t(y,z)x,y,z)$ is a smooth isotopy from $[f_\beta,f_\gamma]$ to the identity, because each $c_t(y,z)$ acts linearly and orthogonally on $S^p$. Schultz's proof of Theorem~\ref{thm:schultz} uses the hypothesis on $\beta$ only to conclude that $f_\beta$ and $f_\gamma$ commute in the group of pseudo-isotopy classes of diffeomorphisms of $S^p\times S^q\times S^r$; the rest of the argument is unchanged.
\end{proof}

In the application we will have $p=6$ and $q=r=3$. The possible Samelson obstruction consequently lies in $\pi_6(\SO(7))$. It is useful to verify the vanishing of this particular group directly, so that the application does not depend on an unstated stability assumption at the edge of the stable range.
\begin{lemma}\label{lem:pi6}
$\pi_6(\SO(7))=0$.
\end{lemma}

\begin{proof}
Consider the fibration $\SO(7)\to\SO(8)\to S^7$. The map $\SO(8)\to S^7$ sends an orthogonal transformation to its value at the unit vector $1$. A section would force the connecting homomorphism in the homotopy sequence to vanish. The normed octonion multiplication gives just such a section. Left multiplication by unit octonions gives a map $S^7\to\SO(8)$, $u\mapsto L_u$, with $L_u(1)=u$. Each $L_u$ preserves the norm, and its determinant is $+1$ by continuity from $u=1$. Thus this is a section of $\SO(8)\to S^7$. Therefore $\pi_7(\SO(8))\to\pi_7(S^7)$ is surjective, the connecting homomorphism $\pi_7(S^7)\to\pi_6(\SO(7))$ vanishes, and $\pi_6(\SO(7))\to\pi_6(\SO(8))$ is injective. Since $6\le 8-2$, the group $\pi_6(\SO(8))$ is in the stable range and equals $\pi_6(\SO)=0$ by Bott periodicity. (See also \cite{Ke}.)
\end{proof}

Schultz's formula gives a homotopy-sphere operation, but does not yet tell us whether its value is zero. To evaluate it we use the description of $\tau$ in smoothing theory. For $n\ge5$ we identify $\Th_n$ with $\pi_n(\mathrm{PL}/\Oo)$.
\begin{lemma}\label{lem:normal-invariant}
Let $\gamma\in\pi_r(\SO(n))$ and $\Sigma\in\Th_n$, $n\ge5$. Then
\begin{equation}\label{eq:normal-invariant}
\chi_{n+r}\bigl(\tau_{r,n}(\gamma,\Sigma)\bigr)=\pm\,\chi_n(\Sigma)\cdot J(\gamma)\in\coker J_{n+r},
\end{equation}
where $J(\gamma)\in\pis r$ is the stable $J$-image and the product is stable composition.
\end{lemma}
\begin{proof}
Schultz describes the Milnor--Munkres--Novikov pairing in smoothing theory as $\tau(\gamma,\Sigma)=\Sigma\circ J(\gamma)$ in $\pi_{n+r}(\mathrm{PL}/\Oo)$, with $J(\gamma)\in\pi_{n+r}(S^n)$ \cite[p.~742]{Sch}. Let $G$ be the monoid of stable self-homotopy equivalences of spheres and $j\colon G\to G/\Oo$ the natural map. Homotopy spheres are stably parallelisable \cite[Theorem~3.1]{KM}, so the image of $\Sigma$ in $\pi_n(G/\Oo)$ lies in $j(\pi_n(G))$, and $\pi_n(G)\cong\pis n$ modulo $\operatorname{im}J_n$ is $\coker J_n\subset\pi_n(G/\Oo)$. Under the Hirsch--Mazur identification $\pi_n(\mathrm{PL}/\Oo)\cong\Gamma_n$, the resulting map $\pi_n(\mathrm{PL}/\Oo)\to\coker J_n$ is the Kervaire--Milnor homomorphism \cite[\S1, pp.~291--292]{Brumfiel}, \cite[\S4]{KM}. Write the image of $\Sigma$ as $j(s)$ with $[s]=\pm\chi_n(\Sigma)$. Precomposition commutes with maps of spaces, so the image of $\Sigma\circ J(\gamma)$ is $j(s)\circ J(\gamma)=j(s\circ J(\gamma))$. Finally $\pi_*(G)\cong\pi_*(Q_0S^0)$ by translation of components, and in $Q_0S^0=\operatorname{colim}_m\Omega^mS^m$ the precomposite $s\circ y$ with $y\in\pi_{n+r}(S^n)$ is adjoint to $S^{n+r+m}\to S^{n+m}\to S^m$; it represents the stable product $s\cdot\Sigma^\infty y$. This proves \eqref{eq:normal-invariant}.
\end{proof}
The same compatibility with stable composition appears in Bredon's work for his pairing $\Th_n\otimes\pis r\to\Th_{n+r}$, $r<n-1$ \cite[Corollary~2.2]{Br}; we do not use it as an input. In dimensions $10$ and $13$ the images of $J$ vanish, so \eqref{eq:normal-invariant} has no indeterminacy.

\subsection{Proof of Theorem B}
We now apply the preceding tools to the explicit construction. First we put the two quotients into the same normalisation of the plumbing pairing. The fixed real axes are useful here: deleting one such axis from each representation gives the reduced orthogonal groups in the definition of $\sigma$.

By Theorem~\ref{thm:sp} applied to $f_{10}=b\ast\mathrm{id}$ and $f_{13}=a_{13}\ast\mathrm{id}$, and since $\mathrm{conj}=\Ad\oplus\mathbf{1}$ and $\R^k\oplus L=(\R^{k-1}\oplus L)\oplus\mathbf{1}$, the two star quotients are plumbing pairings in Schultz's normalisation:
\begin{equation}\label{eq:plumb}
\begin{aligned}
\Sigma^{10}&=\pm\,\sigma_{6,3}(\alpha,\beta),\\
\Sigma^{13}&=\pm\,\sigma_{9,3}(\Ad\circ a_{13},\beta)
=\pm\,\sigma_{9,3}(\alpha\circ g,\beta).
\end{aligned}
\end{equation}
with
\[
\alpha=\Ad\circ b\in\pi_6(\SO(3)),\qquad \beta=\R^2\oplus L\in\pi_3(\SO(6)),\qquad [\beta]=1,
\]
where $\beta$ is stabilised to $\pi_3(\SO(9))$ in the second formula. The equality of the two expressions for $\Sigma^{13}$ follows directly from the actual factorisation $a_{13}=b\circ g$: the first reduced clutching map is $\Ad\circ a_{13}=\alpha\circ g$. The signs in the plumbing formulas only record the orientation conventions.

The factorisation is now in the form required by Schultz: precomposition is by the explicit Hopf construction of $\gamma_0$. We check the range and the commutation condition before applying his formula.

We apply Theorem~\ref{thm:schultz} with $p=6$, $q=3$, $r=3$. Take $\gamma=[\gamma_0]\in\pi_3(\SO(6))$. Lemma~\ref{lem:a13} gives $J(\gamma)=[g]=\nu_6$, hence $\alpha\circ J(\gamma)=[\alpha\circ g]$. The range condition $r=3\le p-2=4$ holds. There is a genuine condition to check before applying the formula. Schultz's original hypothesis fails, since $[\beta]=1$ is not in the image $2\Z$ of $\pi_3(\SO(3))\to\pi_3(\SO(6))$; but the Samelson product of the images of $\beta$ and $\gamma$ in $\pi_3(\SO(7))$ lies in $\pi_6(\SO(7))$, which vanishes by Lemma~\ref{lem:pi6}, so the weakened hypothesis of Theorem~\ref{thm:schultz} holds by Lemma~\ref{lem:samelson}. Hence
\begin{equation}\label{eq:tau}
\Sigma^{13}=\pm\,\tau_{3,10}(\gamma,\Sigma^{10}).
\end{equation}
The same equation also follows from Schultz's original hypothesis alone. The class $2\beta=[\Ad\oplus\R^3]$ lies in the image of $\pi_3(\SO(3))\to\pi_3(\SO(6))$, so Theorem~\ref{thm:schultz} applies to $(\gamma,\alpha,2\beta)$. Bilinearity of $\sigma$ and of $\tau$ then gives $2\,\sigma_{9,3}(\alpha\circ g,\beta)=\sigma_{9,3}(\alpha\circ J(\gamma),2\beta)=\tau_{3,10}(\gamma,\sigma_{6,3}(\alpha,2\beta))=2\,\tau_{3,10}(\gamma,\sigma_{6,3}(\alpha,\beta))$, and multiplication by $2$ is invertible on $\Th_{13}\cong\Z/3$.

It remains to show that the operation in \eqref{eq:tau} is nonzero. By the classical identifications explained in Section~\ref{sec:preliminaries}, $\chi_{10}$ and $\chi_{13}$ are isomorphisms onto $\pi_{10}^s$ and $\pi_{13}^s$, respectively. By Theorem~\ref{thm:sp10}, $\chi_{10}(\Sigma^{10})=\pm\beta_1$. The stable image of $\gamma$ has $J$-image $\nu$. Lemma~\ref{lem:normal-invariant} with $n=10$, $r=3$, together with \eqref{eq:tau}, therefore gives
\[
[\Sigma^{13}]=\pm\,\nu\cdot\beta_1 .
\]
The multiplication can now be checked on its prime-power components. The $3$-primary component of $\nu\in\pis3\cong\Z/24$ is $\alpha_1$, and the $2$-primary component annihilates the $3$-torsion element $\beta_1$. Indeed their product is annihilated by a power of two and by three, and hence is zero by coprimality. Consequently $\nu\cdot\beta_1=\alpha_1\beta_1$, which generates $\pis{13}\cong\Z/3$ \cite[Theorem~1.1.14]{Ra}. Therefore $[\Sigma^{13}]=\pm\alpha_1\beta_1\ne0$, and $\Sigma^{13}$ is exotic. \qed

The calculation separates two kinds of information about the map. It is worth recording this separation because the geometric and topological arguments use different properties of the same formula.
\begin{remark}
The maps used here fit together at both the geometric and topological levels: $a_{13}=b\circ g$ is an actual equivariant factorisation, and $[g]=J(\gamma)$ supplies the ordinary class needed for Schultz's formula. The Hopf construction $g$ represents a generator $\nu_6$; it must not be confused with the threefold suspension of the Blakers--Massey class $\nu'=[b]\in\pi_6(S^3)$, which is $2\nu_6$ after a compatible choice of generators, since $\Sigma\nu'=2\nu_4$. These two choices give the same nonzero stable $13$-class up to sign, but their intermediate coefficients are different \cite{Toda}.
\end{remark}

The role of the six-dimensional trivial summand is now visible. A change of representation may look harmless geometrically, but the second clutching class enters the pairing linearly. The next observation explains why its value was kept equal to a generator.
\begin{remark}\label{rem:choice}
The first clutching class belongs to $\pi_9(\SO(3))\cong\Z/3$ and is stably trivial, since $\pi_9(\SO)\cong\Z/2$. If another equivariant model uses a representation on $U$ whose second reduced clutching class is $3$ times a generator, its plumbing is standard: bilinearity gives $\sigma_{9,3}(u,3\beta)=\sigma_{9,3}(3u,\beta)=0$. This is why retaining the representation $\R^6\oplus L$ matters in the present construction. Existence of a new equivariant map must be checked separately after any representation change.
\end{remark}

\section{The positively curved metric}\label{sec:geometry}

We have proved that the quotient in Section~\ref{sec:construction} has the desired exotic smooth structure. We now construct a metric on that same quotient. The topology determines the attaching map of its two disks, so the metric problem is to produce positively curved disk metrics whose boundary data match under that attaching map.

There are two different local tasks. The southern disk carries a connection metric; its proof must control the interaction between the connection curvature and arbitrary mixed two-planes. The northern disk carries a doubly warped metric; its radius functions are chosen to meet the boundary data already supplied by the south. After comparing their outward second fundamental forms, we smooth the common boundary. We prove all three parts in this section.

The compatible-disk strategy and some of the profiles were motivated by \cite{HLY,GG}. In particular, the transfer statement below has the same scope as \cite[Theorem~4.1]{GG}; our purpose here is to supply the curvature estimates and strict gluing argument used for this explicit $13$-dimensional model. Those existence and transfer theorems are not used as inputs. The ingredients in the proof are Koszul's formula, the Gauss and Codazzi equations, and O'Neill's submersion formula \cite{ON}. We use the curvature and boundary sign conventions of Section~\ref{subsec:curvature-conventions}.

\subsection{Polar trivialisations from invariant parallel transport}
The first task is to choose bundle coordinates which are adapted to the two disks. A general transition function could depend on both the radial and angular variables. That dependence would introduce unnecessary terms in the metric near the boundary. We therefore choose the coordinates by parallel transport along meridians, so that the transition depends only on the angular variable.

Write $S^n=S(\R e\oplus V)$ with the orthogonal $S^3$-action fixing $e$. We use the star-bundle structure of Definition~\ref{def:special-star}. In geodesic polar coordinates, a base point away from the poles is $(\cos t)e+(\sin t)y$, with $y\in S(V)$.

The following lemma supplies the same simple fibre action in both charts. This will later allow the source boundary metrics, the quotient boundary metrics and their second fundamental forms to be compared in one common trivialisation.
\begin{lemma}\label{lem:polar-independent}
The star bundle admits north and south principal trivialisations in which the star action is
\[
q\star(t,y,u)=(t,\rho(q)y,qu),
\]
and the transition has the form
\[
u_S=\theta(y)u_N,\qquad \theta(\rho(q)y)=q\theta(y)q^{-1}.
\]
In particular the transition is independent of $t$.
\end{lemma}
\begin{proof}
We first make parallel transport compatible with the star action. Average any principal connection over the compact star group. The result remains a principal connection because the actions commute. We next choose the fibre coordinate at each pole. At either fixed pole, the star action on the fibre is free. Choose a point in the fibre and identify the fibre with $S^3$ through the principal action. The star action is then given by an injective homomorphism $S^3\to S^3$. It is an automorphism, since its differential is an isomorphism and the groups are connected and simply connected. Every automorphism of $S^3$ is inner, so a choice of the fibre point makes this homomorphism the identity.

Having normalised the action at the pole, we extend that normalisation to the disk. Transport this point along the round radial geodesics from each pole. The resulting section is smooth on the sphere minus the opposite pole. Invariance of the connection and uniqueness of parallel transport give the displayed local star action. On a common meridian both sections are parallel, so their transition is constant along that meridian. Comparing their star actions gives the conjugation identity for $\theta$.
\end{proof}

These coordinates also explain why the quotient pieces are disks. Each star orbit meets the slice with fibre coordinate $u=1$ in exactly one point: apply the star element $u^{-1}$. In Cartesian disk coordinates the resulting invariant map is smooth even at the centre. More explicitly, the quotient of either trivialised polar disk is a standard disk: in Cartesian disk coordinates the quotient map is $(v,u)\mapsto\rho(u)^{-1}v$. Hence comparing the two quotient disk metrics under a common source trivialisation compares them under the actual attaching map, not an arbitrarily chosen boundary marking. Denote the full polar bundle by $P_\theta$.

We also record the attaching map of the two quotient disks, because equality of boundary metrics must be checked under this particular map. For a cut $0<a<\pi$, let $P_N$ and $P_S$ be the restrictions to $\{t\le a\}$ and $\{t\ge a\}$, and put $D_N=P_N/S^3_\star$ and $D_S=P_S/S^3_\star$. On the common boundary, a point with northern coordinates $(y,1)$ has southern fibre coordinate $\theta(y)$. Applying the southern quotient marking therefore gives
\begin{equation}\label{eq:polar-attaching}
\sigma(y)=\rho(\theta(y))^{-1}y,\qquad y\in S(V).
\end{equation}
Conjugation-equivariance of $\theta$ shows that the inverse is $y\mapsto\rho(\theta(y))y$, so $\sigma$ is a diffeomorphism. Thus the full quotient is $D_N\cup_\sigma D_S$. We will write $g_N,g_S$ for the metrics on these disks and $h_N,h_S$ for their induced boundary metrics. The value of $a$ will later be chosen larger than $\pi/2$.

The infinitesimal action records how much an angular direction moves when the group is changed. It is the only representation-dependent quantity in the northern curvature calculation. To measure it, for $y\in S(V)$ put
\[
K_y\xi=\left.\frac d{ds}\right|_{s=0}\rho(e^{s\xi})y,
\qquad \xi\in\im\HH,
\]
with the Euclidean norm on $\im\HH$ and the unit round norm on $S(V)$.

We can now state the geometric assertion that will be applied to our bundle. The number two is a convenient bound for the two quaternionic representations used here. The conclusion concerns the quotient, so in proving it we only need positive curvature on the source two-planes perpendicular to the star orbits.
\begin{theorem}\label{thm:transfer-independent}
If $n\ge3$ and $\|K_y\|_{\mathrm{op}}\le2$ for all $y\in S(V)$, the star quotient $P_\theta/S^3_\star$ admits a smooth metric of strictly positive sectional curvature.
\end{theorem}
The proof is divided into four parts. We first establish a connection-metric estimate valid on every two-plane. We then prove the strict boundary-smoothing statement. After these tools are available we construct the southern and northern disk metrics, compare their boundary data, and complete the proof. This organisation separates the algebraic curvature calculation from the choice of the scalar radius functions.

\subsection{A positive connection metric on a concave base}
We begin with the local metric which will be used on the southern bundle. If the connection were flat and the fibre radius were constant, a mixed plane would have zero curvature. Shrinking such a fibre would not make that plane positively curved. A varying radius provides a new term, $-r^{-1}\operatorname{Hess}r$, in mixed curvature. Our aim is to make this term positive and large enough, relative to the connection errors, while keeping the fibre small.

The Lie algebra metric $Q$ on $\im\HH$ is normalised so that the bi-invariant metric on $S^3$ has sectional curvature one; thus $|[U,V]|=2|U\wedge V|$.

To formulate the estimate we need two bounds on the connection. Its curvature $\Omega$ measures the failure of principal-horizontal directions to commute, and its covariant derivative $D\Omega$ measures the variation of that failure. The derivative uses the base Levi--Civita connection and the induced connection on the adjoint bundle. We use the pointwise norms
\[
|\Omega|^2=\sum_{i<j,a}(\Omega_{ij}^a)^2,
\qquad |D\Omega|^2=\sum_{k,i<j,a}((D_k\Omega)_{ij}^a)^2
\]
in orthonormal frames. Compactness makes their suprema finite.

The parameter $\Lambda$ below is chosen before the small fibre parameter $\varepsilon$. This order is essential. The positive Hessian contribution and one of the mixed error terms are both of order $\varepsilon$, so shrinking $\varepsilon$ alone cannot make their ratio favourable. We first impose enough concavity through $\Lambda$, and only then shrink the fibre to control the remaining errors.
\begin{proposition}\label{prop:connection-independent}
Let $\pi:P\to(B,g_B)$ be a principal right $S^3$-bundle over a compact smooth base, possibly with boundary, and suppose $\sect_{g_B}\ge\kappa>0$. Fix a smooth principal connection $\omega$ with curvature $\Omega$. With the norm convention above, the absolute constant $C=4$ has the following property. Put
\[
M_0=\sup_B|\Omega|,\quad M_1=\sup_B|D\Omega|,
\quad
\Lambda\ge\max\{32CM_0^2,64C^2\kappa^{-1}(M_1+1)^2,1\}.
\]
If $\Phi:B\to\R$ is smooth and $-\operatorname{Hess}\Phi\ge\Lambda g_B$, then there is $\varepsilon_0>0$, depending on the fixed data $\kappa,M_0,M_1,\Lambda,\|\Phi\|_\infty$ and $\|d\Phi\|_\infty$, such that for every $0<\varepsilon<\varepsilon_0$ the metric
\[
G_\varepsilon=\pi^*g_B+r^2Q(\omega,\omega),\qquad
r^2=\varepsilon e^{\varepsilon\Phi},
\]
has strictly positive sectional curvature on every two-plane.
\end{proposition}

\begin{proof}
The proof has three steps. We compute all curvature blocks, combine them using the algebraic relations of an actual two-plane, and finally choose the small parameter. Throughout the proof, horizontal means principal-horizontal, as distinguished from star-horizontal in Section~\ref{subsec:curvature-conventions}.

\emph{The connection and curvature formulas.} Choose a normal base orthonormal frame at the point in question and let $X_i$ be its principal-horizontal lifts. For a fixed $Q$-orthonormal Lie algebra basis $e_a$, put $E_a=r^{-1}e_a^\#$, $\vartheta=d\log r$ and $[e_a,e_b]=c_{ab}^{\ c}e_c$. The fundamental fields correspond to the principal right action. The bracket identities are
\[
[X_i,X_j]^\mathrm{vert}=-r\Omega_{ij}^aE_a,
\qquad [X_i,E_a]=-\vartheta_iE_a,
\qquad [E_a,E_b]=r^{-1}c_{ab}^{\ c}E_c.
\]
The components of $\Omega$ are not constant along a fibre in a fixed Lie algebra basis. Moving along the fibre conjugates the curvature value. Keeping this derivative is necessary for the terms involving the noncommutative Lie bracket of $S^3$. Equivariance of the curvature form gives the identity
\begin{equation}\label{eq:vertical-omega}
E_a(\Omega_{ij}^b)=-r^{-1}c_{ac}^{\ b}\Omega_{ij}^c.
\end{equation}
The bracket identities determine the Levi--Civita connection through Koszul's formula. The first line below describes horizontal differentiation, while the second records the change along the fibre and the bending caused by the varying radius:
\begin{equation}\label{eq:koszul-table}
\begin{aligned}
\nabla_{X_i}X_j&=(\nabla_i^BX_j)^\mathrm{hor}-\tfrac r2\Omega_{ij}^aE_a,
&\nabla_{X_i}E_a&=\tfrac r2\Omega_{ij}^aX_j,\\
\nabla_{E_a}X_i&=\tfrac r2\Omega_{ij}^aX_j+\vartheta_iE_a,
&\nabla_{E_a}E_b&=\tfrac1{2r}c_{ab}^{\ c}E_c-\delta_{ab}\vartheta_iX_i.
\end{aligned}
\end{equation}
To display the positive contribution of the radius separately, write $N=-r^{-1}\operatorname{Hess}r$. This symmetric bilinear form is the mixed curvature term that is already present for a flat connection. The remaining terms below are the corrections caused by $\Omega$, its derivative and the fibre Lie bracket. Differentiating \eqref{eq:koszul-table}, and using \eqref{eq:vertical-omega}, gives the four mixed curvature blocks:
\begin{align}
\langle R(X_i,E_a)E_b,X_j\rangle
&=N_{ij}\delta_{ab}+\tfrac{r^2}{4}\sum_k\Omega_{ik}^b\Omega_{jk}^a
-\tfrac14\langle[e_a,e_b],\Omega_{ij}\rangle,
\label{eq:block-hv}\\
\langle R(X_i,X_j)E_a,E_b\rangle
&=\tfrac{r^2}{4}\sum_k(\Omega_{ik}^a\Omega_{jk}^b-\Omega_{jk}^a\Omega_{ik}^b)
+\tfrac12\langle[e_a,e_b],\Omega_{ij}\rangle,
\label{eq:block-hhvv}\\
\langle R(X_i,E_a)E_b,E_c\rangle
&=\tfrac r2\sum_k\vartheta_k(\delta_{ab}\Omega_{ik}^c-\delta_{ac}\Omega_{ik}^b),
\label{eq:block-hvvv}\\
\langle R(X_i,X_j)X_k,E_a\rangle
&=-\tfrac r2\bigl((D_i\Omega)_{jk}^a-(D_j\Omega)_{ik}^a\bigr)\nonumber\\
&\quad-\tfrac r2(\vartheta_i\Omega_{jk}^a-\vartheta_j\Omega_{ik}^a)
+r\vartheta_k\Omega_{ij}^a.
\label{eq:block-hhhv}
\end{align}
For clarity, the horizontal part of the first calculation, before combining terms, is
\[
\begin{aligned}
\langle\nabla_{X_i}\nabla_{E_a}E_b,X_j\rangle
 &=\tfrac14c_{ab}^{\ c}\Omega_{ij}^c-\delta_{ab}\partial_i\vartheta_j,\\
\langle\nabla_{E_a}\nabla_{X_i}E_b,X_j\rangle
 &=-\tfrac12c_{ac}^{\ b}\Omega_{ij}^c
   +\tfrac{r^2}{4}\sum_k\Omega_{ik}^b\Omega_{kj}^a,\\
\langle\nabla_{[X_i,E_a]}E_b,X_j\rangle
 &=\delta_{ab}\vartheta_i\vartheta_j.
\end{aligned}
\]
These identities give \eqref{eq:block-hv}; the first Bianchi identity then gives \eqref{eq:block-hhvv}. In \eqref{eq:block-hvvv}, the derivative of $c_{ab}^{\ c}/(2r)$ cancels the vertical bracket contribution. In \eqref{eq:block-hhhv}, the last term is the vertical component of $-\nabla_{[X_i,X_j]}X_k$. This explains all derivative and bracket terms in the table.

The four displayed mixed blocks account for the interaction between the two distributions. We also need the two unmixed blocks. The horizontal block is the base sectional curvature with the usual submersion loss. The vertical block is the round fibre curvature with the loss caused by the second fundamental form of the fibre. Explicitly, the purely horizontal sectional numerator and the purely vertical curvature coefficient are, respectively,
\[
\mathcal R(X\wedge Y,X\wedge Y)
=\mathcal R_B(X\wedge Y,X\wedge Y)-\tfrac34r^2|\Omega(X,Y)|^2,
\qquad r^{-2}-|\vartheta|^2.
\]
The first identity also follows from the horizontal submersion formula; the second follows from the Gauss equation and the fibre second fundamental form $-G|_\mathrm{vert}\otimes\nabla\log r$.

\emph{Combining the blocks on an arbitrary two-plane.} Checking the unmixed blocks is not sufficient, because a general tangent two-plane has horizontal and vertical components at the same time. We therefore keep all those components together. Take an arbitrary decomposable bivector
\[
\xi=(X+U)\wedge(Y+V)=\alpha+\beta+\gamma,
\quad
\alpha=X\wedge Y,\quad \beta=X\otimes V-Y\otimes U,\quad\gamma=U\wedge V.
\]
For the off-diagonal contractions we use the symmetric bilinear convention
\[
\mathcal R(A\wedge B,C\wedge D)=\langle R(A,B)D,C\rangle.
\]
Here $X,Y$ are horizontal and $U,V$ vertical. Write $H=|\alpha|$, $M=|\beta|$, and $W=|\gamma|$, so $|\xi|^2=H^2+M^2+W^2$. The three components are orthogonal as bivectors, but their coefficients are not algebraically independent. They come from the same pair of vectors. This dependence is useful: it converts a potentially dangerous linear connection term into a product of the horizontal and vertical areas. The required identity is
\begin{equation}\label{eq:plucker-independent}
\beta_{ia}\beta_{jb}-\beta_{ib}\beta_{ja}=\alpha_{ij}\gamma_{ab}.
\end{equation}
For example, substituting $\beta_{ia}=X_iV_a-Y_iU_a$ makes the terms quadratic in $X$ and quadratic in $Y$ cancel, leaving $(X_iY_j-X_jY_i)(U_aV_b-U_bV_a)$. This is the right-hand side of \eqref{eq:plucker-independent}. With the stated contraction convention, the linear term from \eqref{eq:block-hv} is $-\frac12\langle\Omega(X,Y),[U,V]\rangle$, and the cross term from \eqref{eq:block-hhvv} is $-\langle\Omega(X,Y),[U,V]\rangle$. Their sum is
\[
-\tfrac32\langle\Omega(X,Y),[U,V]\rangle.
\]
Since the norm of the quaternionic bracket is twice the vertical area, its absolute value is at most $3M_0HW$. This is the important cancellation. An error of size $M_0M^2$ would not become small when the fibre shrinks, whereas the positive mixed Hessian term is only of order $\varepsilon$. The actual error couples $H$ with $W$ and can be absorbed using the large vertical curvature.

For completeness we give coefficient bounds for the remaining contractions. They make the estimate independent of an unspecified dimensional constant. Put $\lambda=\lambda_{\min}N$. The purely horizontal and vertical terms and the four mixed blocks give
\begin{align}
\mathcal R_G(\xi,\xi)\ge{}&
(\kappa-\tfrac34r^2M_0^2)H^2+
(\lambda-\tfrac12r^2M_0^2)M^2+(r^{-2}-|\vartheta|^2)W^2\nonumber\\
&-2r(M_1+2|\vartheta|M_0)HM
-(3M_0+r^2M_0^2)HW-2r|\vartheta|M_0MW.
\label{eq:connection-sharp-bound}
\end{align}
Here is the coefficient check. The full array $\Omega_{ij}^a$ has norm $\sqrt2|\Omega|$. Hence the quadratic array $\sum_k\Omega_{ik}^b\Omega_{jk}^a$ in \eqref{eq:block-hv} has norm at most $2|\Omega|^2$, and its contraction contributes at worst $r^2M_0^2M^2/2$. The quadratic part in \eqref{eq:block-hhvv} has full-array norm at most $r^2|\Omega|^2$. It is antisymmetric in $i,j$ and in $a,b$, so restriction to $i<j$, $a<b$ divides this norm by two; the cross-block factor two restores the bound $r^2M_0^2HW$. Its bracket term has already been combined with the mixed diagonal term above, giving $3M_0HW$.

The full array in \eqref{eq:block-hhhv} has norm at most $\sqrt2rM_1+2\sqrt2r|\vartheta|M_0$. Restricting to $i<j$ divides by $\sqrt2$, and the cross-block factor two gives the $HM$ term in \eqref{eq:connection-sharp-bound}. Finally put $w_i^a=\sum_k\vartheta_k\Omega_{ik}^a$. Then
\[
\sum_{i,a}|w_i^a|^2\le2|\vartheta|^2|\Omega|^2,\qquad
\sum_{a,b,c}(\delta_{ab}w_i^c-\delta_{ac}w_i^b)^2=4|w_i|^2.
\]
These inequalities give full-array norm at most $\sqrt2r|\vartheta|M_0$ for \eqref{eq:block-hvvv}. Restriction to $b<c$ and the cross-block factor two give the final term $2r|\vartheta|M_0MW$. This proves \eqref{eq:connection-sharp-bound} by Cauchy--Schwarz. Only the sectional lower bound on the base was used, since $\alpha=X\wedge Y$ is decomposable.

For the parameter choices it is convenient to use one common coefficient. Each term in \eqref{eq:connection-sharp-bound} is bounded by the corresponding term below with $C=4$:
\begin{align}
\mathcal R_G(\xi,\xi)\ge{}&
(\kappa-Cr^2M_0^2)H^2+
(\lambda-Cr^2M_0^2)M^2+(r^{-2}-|\vartheta|^2)W^2\nonumber\\
&-Cr(M_1+|\vartheta|M_0)HM
-C(M_0+r^2M_0^2)HW-Cr|\vartheta|M_0MW.
\label{eq:connection-bound}
\end{align}
Thus the common coefficient is explicit; the less sharp form will simplify the successive Young inequalities.

\emph{The choice of the radius and the small parameter.} We now use the special form of $r$ in the proposition. Its logarithmic derivative is of order $\varepsilon$, while its vertical curvature is of order $\varepsilon^{-1}$. These different orders are what allow the estimates to close. Put $D_0=\sup_B|d\Phi|$. For our radius function,
\[
\vartheta=\tfrac\varepsilon2d\Phi,
\qquad
N=-\tfrac\varepsilon2\operatorname{Hess}\Phi
-\tfrac{\varepsilon^2}{4}d\Phi\otimes d\Phi.
\]
First choose $\varepsilon$ so small that $|\varepsilon\Phi|\le\log2$, $|\vartheta|M_0\le1$, $2\varepsilon M_0^2\le1$, and
\[
4C\varepsilon M_0^2\le\kappa,
\qquad \varepsilon D_0^2\le\Lambda/4,
\qquad \varepsilon^3D_0^2\le1.
\]
The first smallness condition protects the positive base curvature. The next two control the gradient-square losses in the Hessian and vertical terms. Together with the lower bound chosen for $\Lambda$, they imply that the three diagonal coefficients in \eqref{eq:connection-bound} are at least $\kappa/2$, $3\varepsilon\Lambda/8$ and $1/(4\varepsilon)$. The first error is bounded by
\[
Cr(M_1+|\vartheta|M_0)HM
\le\tfrac\kappa8H^2+4C^2\varepsilon\kappa^{-1}(M_1+1)^2M^2
\le\tfrac\kappa8H^2+\tfrac{\varepsilon\Lambda}{16}M^2.
\]
The other errors have a different balance. The horizontal--vertical area error is absorbed by the large $W^2$ coefficient, while the last error is of higher order because $\vartheta=O(\varepsilon)$. For these second and third errors Young's inequality gives, on decreasing $\varepsilon$ further,
\[
\begin{aligned}
C(M_0+r^2M_0^2)HW
&\le\tfrac\kappa8H^2+\tfrac1{16\varepsilon}W^2,\\
Cr|\vartheta|M_0MW
&\le\tfrac1{16\varepsilon}W^2+\tfrac{\varepsilon\Lambda}{16}M^2.
\end{aligned}
\]
Indeed it suffices additionally to impose
\[
\varepsilon\le\frac{\kappa}{32C^2(M_0+1)^2},
\qquad32C^2\varepsilon^3D_0^2M_0^2\le\Lambda.
\]
Every condition is an upper bound on $\varepsilon$ after all other data have been fixed. Thus a common positive choice exists, and
\begin{equation}\label{eq:south-uniform}
\mathcal R_{G_\varepsilon}(\xi,\xi)
\ge\tfrac\kappa4H^2+\tfrac{\varepsilon\Lambda}{4}M^2+\tfrac1{8\varepsilon}W^2.
\end{equation}
A nonzero two-plane has at least one of $H,M,W$ nonzero. Since all three coefficients are positive for the fixed admissible $\varepsilon$, \eqref{eq:south-uniform} is positive on every such plane. In fact it gives a uniform sectional-curvature lower bound by the minimum of those three coefficients. This bound is for that fixed metric; it is not asserted to stay positive uniformly as $\varepsilon\to0$. This is the all-plane conclusion needed for the southern metric. It proves the proposition.
\end{proof}

\subsection{An elementary strict gluing lemma}
Before constructing the two pieces, we explain the boundary condition that we will arrange. Equality of boundary metrics gives a continuous metric on their union. The first normal derivatives may still differ, so this joined metric need not be smooth. The sum of the outward second fundamental forms measures that difference. When it is strictly positive, the jump can be spread over a short collar so as to create a large positive radial curvature term.

The following lemma proves the precise strict case needed in this paper. It is also a special case of the positive-sectional-curvature conclusion of \cite[Theorem~A(i), $k=1$]{RW}, which permits a semidefinite sum. We give the strict-case proof so that the present construction does not depend on that more general smoothing result.
\begin{lemma}\label{lem:strict-glue-independent}
Let two compact Riemannian manifolds have strictly positive sectional curvature, isometric boundaries, and outward second fundamental forms $B_-$ and $B_+$ whose sum, under that isometry, is strictly positive definite. Their union admits a smooth metric of strictly positive sectional curvature, equal to the given metrics outside an arbitrarily small collar of the seam.
\end{lemma}
\begin{proof}
We first choose one normal coordinate across the seam. It increases towards the boundary on the first side and into the interior on the second. This accounts for the opposite signs of the two outward second fundamental forms. Identify the boundaries as a compact manifold $S$. Write the geodesic collar metrics as $dt^2+h_-(t)$ for $t\le0$ and $dt^2+h_+(t)$ for $t\ge0$, with common value $h_0$. Put
\[
a=h_-'(0)=2B_-,\qquad b=h_+'(0)=-2B_+,\qquad D=a-b\ge d h_0>0.
\]
The interpolation must match both the metric and its first derivative at its endpoints. A cubic polynomial has exactly enough coefficients for these four conditions. On $[-\eta,\eta]$ use the tensor-valued cubic Hermite polynomial $H_\eta(t)$ matching the values and first derivatives at $-\eta$ and $\eta$. Outside it keep the original metrics. For the linear germs $h_0+ta$ and $h_0+tb$ this polynomial is exactly
\[
h_0+\tfrac12(a+b)t-\tfrac{t^2+\eta^2}{4\eta}D.
\]
Taylor expansion of the actual germs therefore gives, with $\tau=(t+\eta)/(2\eta)$,
\begin{equation}\label{eq:hermite-jets}
H_\eta=h_0+O(\eta),\quad
H_\eta'=(1-\tau)a+\tau b+O(\eta),\quad
H_\eta''=-\frac{D}{2\eta}+O(1).
\end{equation}
These estimates hold with two, one and zero spatial derivatives, respectively; all constants are uniform on $S$. The resulting metric is positive definite, $C^1$ and piecewise smooth.

The metric is now $C^1$, but that alone does not imply positive curvature. We check the tangential, radial and cross contributions before performing the final smoothing. The strict boundary inequality is used in two ways: it makes the determinant comparison favourable for tangential planes, and it makes radial curvature large and positive. For an $h_0$-orthonormal tangent pair $X,Y$, restrict $a/2$ and $b/2$ to their two-plane and call the resulting matrices $A,B$. The Gauss equation and \eqref{eq:hermite-jets} give the tangential numerator
\[
\mathcal R_{h_0}(X\wedge Y,X\wedge Y)
-\det((1-\tau)A+\tau B)+O(\eta).
\]
For tangential planes, the intrinsic boundary curvature stays close to that of $h_0$. The term that needs control is the determinant of the interpolated second fundamental form. The following elementary two-by-two identity makes the sign visible:
\[
\det((1-\tau)A+\tau B)
=(1-\tau)\det A+\tau\det B-\tau(1-\tau)\det(A-B)
\]
and $A-B>0$ show that this numerator is bounded below by the convex combination of the original two tangential numerators, minus $O(\eta)$. Compactness and the original strict positive sectional curvature make it at least $c>0$ for small $\eta$, uniformly on tangent two-planes.

For a plane containing the normal direction, the dominant term is instead the second normal derivative of the metric. The cubic interpolation has made this derivative equal to $-D/(2\eta)$ to leading order. The exact radial identity is
\[
\langle R(\partial_t,Z)Z,\partial_t\rangle
=-\tfrac12H_\eta''(Z,Z)
+\tfrac14(H_\eta'H_\eta^{-1}H_\eta')(Z,Z)
\ge\frac d{8\eta}|Z|_{h_0}^2
\]
for small $\eta$. The Codazzi block contains only spatial derivatives of $H_\eta'$ and lower jets, so it is bounded uniformly as $\eta\to0$. We now combine these estimates, rather than restricting the proof to planes either tangent or normal to the slices. Every decomposable bivector has the form
\[
\xi=(\lambda\partial_t+X)\wedge(\mu\partial_t+Y)
=\partial_t\wedge Z+\alpha,\quad Z=\lambda Y-\mu X,\quad\alpha=X\wedge Y.
\]
Its tangential component is still decomposable. Thus, for a uniform constant $C_0$,
\[
\mathcal R(\xi,\xi)\ge c|\alpha|_{h_0}^2+\frac d{8\eta}|Z|_{h_0}^2
-2C_0|\alpha|_{h_0}|Z|_{h_0}>0
\]
for small enough $\eta$, by Young's inequality. Fix such an $\eta$. This also gives a uniform positive lower bound in a fixed background norm on both sides of the two $C^1$ interfaces at $\pm\eta$.

The preceding estimate gives strict positive curvature on each smooth piece and on both one-sided limits at the two interfaces. We must still remove the jumps in second derivatives. A general smoothing argument would require care because curvature is nonlinear in the metric. Here its dependence on the only discontinuous derivative is affine, which is enough to control the smoothing. Convolve the piecewise smooth tensor field $H$ in $t$ with a nonnegative smooth mollifier of scale $\delta$, leaving the $S$ variable unchanged. Its lower jets $H,H',\nabla_SH,\nabla_S^2H,\nabla_SH'$ converge uniformly: they are continuous because the endpoint values and first derivatives matched as smooth tensor fields on $S$. Since $H'$ is continuous, the distributional $H''$ has no Dirac mass; it is a bounded piecewise smooth field for this fixed $\eta$. The covariant curvature tensor of $dt^2+H$ depends affinely on $H''$, whose only occurrence is $-H''/2$ in the radial block, and continuously on the other displayed jets. Consequently, in the fixed splitting $\R\partial_t\oplus TS$,
\[
R(dt^2+H*\rho_\delta)(t)
=\int\rho_\delta(s)R(dt^2+H)(t-s)\,ds+o(1)
\]
uniformly as $\delta\to0$, with $\eta$ fixed. More explicitly, use $dt^2+h_0$ as the fixed background metric and let $c_\eta>0$ be the lower bound on its unit decomposable bivectors obtained after fixing $\eta$. Choose $\delta$ so that the displayed tensor error has operator norm on bivectors less than $c_\eta/2$. For each fixed decomposable bivector, positivity of its numerator is a linear inequality in the covariant curvature tensor. Taking all these inequalities together gives a convex cone. Averaging tensors already satisfying the same strict lower bound therefore preserves that bound, and the smoothed numerator is at least $c_\eta/2$ on these unit bivectors. No estimate uniform in a simultaneous limit of $\eta$ and $\delta$ is needed.

We finally make sure that the smoothing changes only a small neighbourhood of the interfaces. This matters because the two disk metrics should be retained away from the gluing collar. For a local modification, multiply $(H*\rho_\delta)-H$ by cutoffs supported in disjoint small neighbourhoods of $\pm\eta$ and equal to one near each interface. Fix these neighbourhoods and cutoffs first, and choose the convolution scale smaller than the distance from their transition regions to the interfaces. On those transition regions the original $H$ is smooth, and convolution converges in $C^2$; all cutoff derivative errors consequently tend to zero there. Positivity is open in $C^2$, so it is preserved on these regions as well. The resulting metric is smooth and positive everywhere and unchanged outside the chosen collar. Here the union carries the smooth structure given by the two normal collars; by uniqueness of collars it is diffeomorphic, by a diffeomorphism which is the identity away from the seam, to the union with any other compatible smooth structure. This proves the strict gluing statement directly.
\end{proof}

\subsection{The southern disk}
We now build the southern metric and record the boundary data that the north must match. We take the southern base to be a round cap of radius less than $\pi/2$. On such a cap a round height function is strictly concave, which is exactly the hypothesis needed in Proposition~\ref{prop:connection-independent}.

Choose fixed numbers $\pi/2<a<b_0<\pi$ and a smooth $\chi:[a,\pi]\to[0,1]$ equal to zero near $a$ and one for $t\ge b_0$. The connection has to satisfy two local requirements. Near the common boundary it should be the product connection in the northern trivialisation, so that the northern metric can match it. Near the south pole it should be the product connection in the southern trivialisation, so that it is smooth at the centre. The cutoff interpolates between these requirements while respecting the transition function. On the southern round cap $U_S=\{t\ge a\}$ use the connection with local forms
\[
A_N=\chi\theta^{-1}d\theta,\qquad
A_S=(\chi-1)d\theta\,\theta^{-1},\qquad
\omega=u^{-1}Au+u^{-1}du.
\]
The coordinate change $u_S=\theta u_N$ requires $A_N=\theta^{-1}A_S\theta+\theta^{-1}d\theta$, which the formulas satisfy. Near the south pole $A_S=0$, so the connection is smooth. Near the boundary $A_N=0$. To check star-invariance explicitly, conjugation-equivariance of $\theta$ gives $\rho(q)^*A_N=qA_Nq^{-1}$. Since the star action also sends $u$ to $qu$, its pullback on the full connection form is
\[
(qu)^{-1}(qA_Nq^{-1})(qu)+(qu)^{-1}d(qu)
=u^{-1}A_Nu+u^{-1}du=\omega.
\]
The calculation in the southern chart is identical. Its curvature and covariant derivative therefore have finite bounds $M_0,M_1$ on the compact cap.

The two curvature bounds of this fixed connection are now determined. We can therefore choose the concavity before choosing the fibre scale, as required by the proposition. Choose $\Lambda$ as in Proposition~\ref{prop:connection-independent} with $\kappa=1$, choose $A_0>0$ so that $A_0|\cos a|\ge\Lambda$, and set
\[
\Phi(t)=-A_0\cos t,\qquad
G_S=dt^2+\sin^2t\,h_{S^{n-1}}+r_S(t)^2Q(\omega,\omega),
\quad r_S^2=\varepsilon e^{-\varepsilon A_0\cos t}.
\]
Since $\operatorname{Hess}(\cos t)=-\cos t\,g_B$, we have $-\operatorname{Hess}\Phi=-A_0\cos t\,g_B\ge\Lambda g_B$. The proposition proves $\sect_{G_S}>0$ for small $\varepsilon$. At the centre put $\zeta=\pi-t$; $r_S=\sqrt\varepsilon\exp(\varepsilon A_0\cos\zeta/2)$ is smooth and even, and the connection is a product there. This verifies smoothness at the centre. The star action is free and isometric. O'Neill's formula therefore gives $\sect_{g_S}>0$ for the quotient disk metric, including centre and boundary.

For the next part of the proof, the important output is not an explicit formula for the whole quotient metric. It is the boundary metric and its first normal derivative. The product form of the connection near the boundary makes both quantities explicit. The source boundary in the common northern trivialisation is the product metric
\[
F_a^2h_{S^{n-1}}+r_a^2h_{S^3},\quad
F_a=\sin a,\quad r_a^2=\varepsilon e^{-\varepsilon A_0\cos a}.
\]
The logarithmic derivatives describe the two principal contributions to the second fundamental form: one from the angular sphere and one from the fibre. We denote them by
\[
q_s=\frac{r_S'(a)}{r_a}=\frac\varepsilon2 A_0F_a,\qquad
\mu_S=\cot a<0.
\]
For scaled orthonormal angular and fibre components $X,U$, the outward source second fundamental form is
\begin{equation}\label{eq:south-II}
B_S^{\mathrm{src}}(X+U,X+U)=-\mu_S|X|^2-q_s|U|^2.
\end{equation}

\subsection{The northern disk and all its horizontal two-planes}
The northern metric is chosen to match the southern boundary data, which at this stage form a family depending on $\varepsilon$. We keep $\varepsilon$ variable while choosing $F$ and $\ell$, and then define the fibre profile for each sufficiently small $\varepsilon$. It has two radius functions: $F$ controls the angular sphere of the disk and $r$ controls the principal fibre. The first must collapse smoothly at the disk centre; the second must stay positive there because the fibre does not collapse in the total space. At the other end both radii must match the southern values. We will take the northern fibre logarithmic derivative to be twice the southern one, giving the strict boundary inequality an additional positive term.

Identify the northern disk with radial coordinate $s$ by the equivariant radial diffeomorphism $t=(a/\ell)s$, preserving the angular variable; this is linear at the centre. On its trivial bundle use
\[
G_N=ds^2+F(s)^2h_{S^{n-1}}+r(s)^2h_{S^3},\qquad0\le s\le\ell.
\]
The total-space metric will be checked on star-horizontal planes, which are the planes seen by the quotient. The two product factors are not separately perpendicular to the star action, so we first determine the actual horizontal space. At fibre coordinate $u=1$, the star orbit vector generated by $\xi$ is $(FK_y\xi,r\xi)$ in scaled components. A star-horizontal vector is therefore
\[
\lambda\partial_s+X+U,\qquad U=TX,\qquad T=-\frac F rK_y^*,\qquad\|T\|\le\frac{2F}r.
\]
Principal right translation is an isometry commuting with the star action, so this describes every fibre point. The graph formula is useful because the fibre component of a horizontal vector is controlled by its angular component. It therefore allows negative curvature terms involving the fibre to be estimated by positive terms involving the base. For a pair of such horizontal vectors with components $(\lambda,X,U)$ and $(\mu,Y,V)$,
\begin{equation}\label{eq:graph-estimates-independent}
|\lambda V-\mu U|^2\le\frac{4F^2}{r^2}|\lambda Y-\mu X|^2,
\qquad
|X\otimes V-Y\otimes U|^2\le\frac{8F^2}{r^2}|X\wedge Y|^2.
\end{equation}
The first estimate is immediate from $U=TX,V=TY$. For the second, a determinant-one change of the pair makes $X,Y$ orthogonal while preserving both alternating expressions. Their squared tensor norm is then $|X|^2|TY|^2+|Y|^2|TX|^2$, proving the estimate; dependent pairs give zero on both sides.

We next compute the numerator on a completely arbitrary horizontal plane. Since the two factors are round and their radii depend only on $s$, the curvature tensor has five types of diagonal blocks, corresponding to radial--angular, radial--fibre, angular, fibre, and angular--fibre bivectors. In radially parallel orthonormal frames, their coefficients are
\[
-\frac{F''}F,\quad-\frac{r''}r,\quad
\frac{1-F'^2}{F^2},\quad\frac{1-r'^2}{r^2},\quad-\frac{F'r'}{Fr}.
\]
They follow from Koszul's formula: $\nabla_{E_i}\partial_s=(F'/F)E_i$, $\nabla_{E_a}\partial_s=(r'/r)E_a$, the radial components of $\nabla_{E_i}E_j$ and $\nabla_{E_a}E_b$ are $-(F'/F)\delta_{ij}\partial_s$ and $-(r'/r)\delta_{ab}\partial_s$, and the two factors have intrinsic curvature one. All other off-diagonal curvature blocks vanish. Thus the numerator on the arbitrary horizontal pair is exactly
\begin{align}
\mathcal K={}&-\frac{F''}F|\lambda Y-\mu X|^2
-\frac{r''}r|\lambda V-\mu U|^2\nonumber\\
&+\frac{1-F'^2}{F^2}|X\wedge Y|^2
+\frac{1-r'^2}{r^2}|U\wedge V|^2
-\frac{F'r'}{Fr}|X\otimes V-Y\otimes U|^2.
\label{eq:north-all-planes}
\end{align}
The last formula explains what the scalar profiles have to accomplish. Positive angular and radial--angular curvature must dominate the two terms which may have the wrong sign. On the punctured disk $0<s\le\ell$, the graph estimates show that it is enough to impose the following scalar inequalities. Under these inequalities the numerator is strictly positive on every nonzero horizontal decomposable bivector:
\begin{equation}\label{eq:profile-targets}
F,r>0,\quad F',r'\ge0,\quad r'<1,\quad
\frac{1-F'^2}{F^2}>\frac{8FF'r'}{r^3},\qquad
-\frac{F''}F>\frac{4F^2(r'')_+}{r^3}.
\end{equation}
Indeed \eqref{eq:graph-estimates-independent} absorbs the two negative terms in \eqref{eq:north-all-planes}. If both remaining radial-angular and angular-area terms vanished, $(\lambda,X)$ and $(\mu,Y)$ would be dependent, and hence the original graph vectors would be dependent. This excludes a zero plane. At the centre we instead require $F(0)=0$, $F'(0)=1$, an odd smooth extension of $F$, and $r(0)>0$ with $r$ constant nearby. The resulting smooth metric and its curvature at the centre are verified below; expressions with $F$ in the denominator are only used for $s>0$.

We now choose the profiles. The small parameter $\delta$ controls the concavity of $F$. The differential equation for $F$ gives a smooth collapse at the centre and keeps $F'$ positive all the way to the boundary. For $r$, we choose a boundary logarithmic derivative $q_N=2q_s$, twice the southern value. The extra positive derivative will leave a positive fibre term in the sum of outward second fundamental forms. We switch that derivative off near the centre by a cutoff, so that $r$ is constant there. The transition is placed where $F$ is still small, making the possible radial--fibre loss easy to absorb.

Here are profiles satisfying all conditions with the required boundary data. Choose smooth nondecreasing $\eta:\R\to[0,1]$, zero for $s\le1/4$ and one for $s\ge1/2$, and let $C_\eta=\sup|\eta'|$. Choose $\delta>0$ with
\[
\delta\le F_a/2,\qquad
\delta^3\le\frac1{256A_0F_a(1+C_\eta)}.
\]
The equation for $F$ is chosen so that both its relevant curvature coefficients can be read off explicitly. Define
\[
F'=e^{-F^2/(2\delta^2)},\quad F(0)=0,\qquad
\ell=\int_0^{F_a}e^{z^2/(2\delta^2)}\,dz,
\]
so $F(\ell)=F_a$. Put $q_N=2q_s=\varepsilon A_0F_a$, $d=r_aq_N$ and
\[
r(s)=r_a-d\int_s^\ell\eta(v/\delta)\,dv.
\]
The length $\ell$ may be large, but it is finite and does not depend on $\varepsilon$. No upper bound on this fixed length is needed: we can subsequently take $\varepsilon$ small enough to meet the condition involving $q_N\ell$. The cutoff, $\delta$, $F$ and $\ell$ are now fixed independently of $\varepsilon$. The quantities $r_a,q_s,q_N,d$ and the function $r$ still depend on $\varepsilon$. Since $q_N=O(\varepsilon)$, $d=O(\varepsilon^{3/2})$, decrease $\varepsilon$ so that $q_N\ell\le1/2$ and $d<1$. Then
\[
r_a/2\le r\le r_a,\quad r'=d\eta(s/\delta),\quad
r''=(d/\delta)\eta'(s/\delta),\quad
\frac d{r^3}\le\frac{8q_N}{r_a^2}\le8A_0F_a.
\]
Moreover
\[
-\frac{F''}F=\delta^{-2}e^{-F^2/\delta^2},\qquad
\frac{1-F'^2}{F^2}=\frac{1-e^{-F^2/\delta^2}}{F^2}.
\]
We check the two target inequalities separately. The only region where $r''$ can cause a loss is the short interval in which the cutoff changes. There the bound $F\le s$ keeps the angular radius small, while the concavity of $F$ remains uniformly strong. On the support of $r''$, $F\le s\le\delta/2$, and hence
\[
-\frac{F''}F\ge e^{-1/4}\delta^{-2},\qquad
\frac{4F^2r''}{r^3}\le8A_0F_aC_\eta\delta
\le\frac{C_\eta}{32(1+C_\eta)}\delta^{-2}.
\]
This proves the radial--angular domination where a loss is present. Off that support the second target in \eqref{eq:profile-targets} is immediate. The angular inequality must hold on the whole interval, including where $F$ is not small compared with $\delta$. We use a one-variable estimate which is uniform in their ratio. Elementary calculus gives
\[
\frac{2\sinh(x^2/2)}{x^3}\ge\frac1x+\frac{x^3}{24}
\ge\frac4{3\,8^{1/4}}>\frac12\qquad(x>0).
\]
With $x=F/\delta$, this yields
\[
\frac{1-F'^2}{F^2}>\frac{FF'}{2\delta^3}
>64A_0F_aFF'\ge\frac{8FF'r'}{r^3}.
\]
Also $r'\le d<1$, proving all targets.

The profile inequalities prove positivity away from the centre. We also have to check that the polar formulas describe a smooth metric there, rather than merely a metric on a punctured disk. This is why the oddness of $F$ and the constancy of $r$ near zero were built into the construction. The inverse of the odd smooth function $z\mapsto\int_0^ze^{v^2/(2\delta^2)}dv$ is odd and smooth. It gives
\[
F(s)=s-\frac{s^3}{6\delta^2}+O(s^5).
\]
Since $r$ is constant near zero, $G_N$ extends smoothly across its centre. There the star orbit is the whole fibre, and the star-horizontal tangent space is the base tangent space, with source sectional curvature $\delta^{-2}$. Thus source horizontal two-planes have positive curvature even at the centre. The star action is free and isometric, and O'Neill's formula gives $\sect_{g_N}>0$ on the entire quotient disk. The infinitesimal base action can have smaller rank at special points. This does not cause a singularity in the quotient, because the star action still moves the fibre freely. More quantitatively, the Gram matrix of orbit vectors is $F^2K_y^*K_y+r^2I\ge r^2I$.

\subsection{Boundary comparison and completion of the proof}
Both quotient metrics are now smooth and strictly positively curved. To join them we must use the actual identification coming from the bundle, and verify two conditions under that identification: equality of the induced metrics and positivity of the sum of outward second fundamental forms. The common source trivialisation is what makes these checks direct.

At the northern boundary $F=F_a$, $r=r_a$, $r'/r_a=q_N=2q_s$, and $\mu_N=F'(\ell)/F_a>0$. Hence
\[
B_N^{\mathrm{src}}(X+U,X+U)=\mu_N|X|^2+2q_s|U|^2.
\]
The source boundary metric equals the southern one exactly, in the common northern trivialisation; the star actions are also exactly the same. Consequently their quotient boundary metrics agree under the actual attaching map $\sigma$.

For the second fundamental forms, it is best to use the same horizontal lift on both sides. Their source metrics and star actions already agree at the boundary, so the horizontal spaces there are the same. For a quotient tangent vector $Y$, take its common source star-horizontal lift $X+U$. The normal fields $\partial_s$ and $-\partial_t$ are star-horizontal (and principal-horizontal). The quotient second fundamental form equals the source second fundamental form on horizontal lifts, by the horizontal part of the submersion connection. Using \eqref{eq:south-II}, the choice $q_N=2q_s$ leaves a strictly positive fibre term:
\[
(B_N+\sigma^*B_S)(Y,Y)=(\mu_N-\mu_S)|X|^2+q_s|U|^2.
\]
Both coefficients are positive: $\mu_N>0>\mu_S$ and $q_s=\varepsilon A_0F_a/2>0$. Since $Y$ has squared quotient norm $|X|^2+|U|^2$, this immediately gives the strict uniform bound
\begin{equation}\label{eq:strict-boundary-final}
h_N=\sigma^*h_S=h,\qquad
B_N+\sigma^*B_S\ge
\min\{\mu_N-\mu_S,q_s\}\,h>0.
\end{equation}
This is also the restriction of a positive definite source boundary form, so no separate lower bound for the angular component of a horizontal lift is needed at this step. The graph estimates remain essential for the northern curvature calculation.

It remains to check that all choices made in the proof can be satisfied at once. The order is part of the argument. The connection bounds are fixed before $A_0$, the angular profile $F$ and the length $\ell$ are fixed before the final fibre scale is chosen, and the smoothing scale is chosen only after both metrics are fixed. Explicitly, choose parameters in the order just used: $a,b_0,\chi$, the connection bounds, $\Lambda$ and $A_0$, then $\eta,\delta,\ell$, and finally a positive $\varepsilon$ satisfying the finitely many smallness conditions. This order has no circular dependency. Both quotient disk metrics are now fixed, smooth and strictly positively curved. Apply Lemma~\ref{lem:strict-glue-independent} to \eqref{eq:strict-boundary-final}. Its smoothing scale is chosen after these metrics are fixed. This proves Theorem~\ref{thm:transfer-independent}.

The transfer theorem was proved for any polar star bundle satisfying the stated action-field bound. To finish the main geometric assertion, we only need to verify that the representation used in Section~\ref{sec:construction} has that bound. The explicit form of the representation makes this last verification short.
\begin{proof}[Proof of Theorem A]
For the bundle constructed in Section~\ref{sec:construction}, the two actions commute and the star action is free because it is pulled back equivariantly from the free Gromoll--Meyer star action. Lemma~\ref{lem:polar-independent} supplies its polar presentation. The representation fixes the real axis of the last quaternionic summand. On its orthogonal complement $V=\R^6\oplus\HH\oplus\im\HH$, for $y=(z,w,y')$ and $\xi\in\im\HH$,
\[
K_y\xi=(0,\xi w,[\xi,y']),\qquad
|K_y\xi|^2=|\xi|^2|w|^2+4|\xi\times y'|^2\le4|\xi|^2.
\]
Thus the independently proved transfer theorem applies and gives the required smooth metric of strictly positive sectional curvature.
\end{proof}

\section*{Verification}
Selected finite-dimensional algebraic identities and scalar inequalities used
in the southern and northern curvature estimates and in the gluing argument
have been checked in Lean~4. The accompanying source files specify the
formalised statements and their assumptions. The curvature formulas, smooth
extension and quotient constructions, analytic smoothing, and topological
identification are not formalised in those files. These certificates do not
constitute a formal proof of Theorems~A or~B.

\section*{Acknowledgements}
This work was supported by the National Natural Science Foundation of China
(NSFC), grant no.\ 12131012.

Generative AI tools (ChatGPT/Codex and Claude) were used in mathematical
exploration and checking, drafting and revision, reference checking, and the
preparation of verification code. Responsibility for the mathematical content
rests with the named authors.

\end{document}